\documentclass[11pt,a4paper,leqno]{article}
\usepackage{latexsym}
\usepackage{amsfonts}
\usepackage {amsbsy}
\usepackage {amsmath}
\usepackage{amssymb}
\newtheorem{defn}{Definition}[section]

\newtheorem{lem}[defn]{Lemma}

\newcommand \demo{ Proof: }
\newcommand \C{\mathbb C}

\newcommand \R{\mathbb R}

\newcommand \fin{$\blacktriangleright$\\}
\newcommand \lto{\longmapsto}
\newcommand \lra{\longrightarrow}

\newcommand \sub{\subset}

\newcommand \al{\alpha}
\newcommand \be{\beta}

\newcommand \ep{\varepsilon}
\newcommand \vf{\varphi}
\newcommand \si{\sigma}

\newcommand \om{\omega}

\newcommand \we{\wedge}

 \title{A priori $L^{\infty}$-estimates for degenerate complex 
 Monge-Amp\`ere equations}

 \author{P. Eyssidieux, V. Guedj and A. Zeriahi} 
\begin{document}
\maketitle
  
\noindent{\bf Abstract :} We study families of complex Monge-Amp\`ere equations, focusing on the case where the 
cohomology classes degenerate to a non big class. 
We establish uniform a priori $L^{\infty}$-estimates for the normalized solutions, generalizing the recent work of S. Kolodziej and G. Tian. This has interesting consequences in the study of the K\"ahler-Ricci flow.
\vskip 0.3 cm
 \section{Introduction}
 Let $\pi : X \lra Y$ be a non degenerate holomorphic mapping between compact K\"ahler manifolds such that $n := dim_{\C} X \geq m := dim_{\C}Y$. Let $\omega_X$, $\omega_Y$  K\"ahler forms on $X$ and $Y$ respectively.
 Let $F : X \lra \R^+$ be a non negative function such that $F \in L^p (X)$ for some $p > 1$.

 Set $\omega_t := \pi^* (\omega_Y) + t  \omega_X$, $t>0$.
 We consider the following  family of complex Monge-Amp\`ere equations
 
\begin{displaymath}
 \left \{\begin{array}{ll}
 (\omega_t + dd^c \vf_t)^n & = \ c_t t^{n - m} F \omega_X^n  \\
  \max_X \vf_t = 0 & =  \ 0  \\
 \end{array} \right. \leqno (\star)_t
 \end{displaymath}
 where $\vf_t$
 is $\omega_t-$plurisubharmonic on $X$ and $c_t > 0$ is a constant given by
 $$ c_t t^{n - m} \int_X F \omega_X^n = \int_X \omega_t^n.$$
 
  It follows from the seminal work of S.T. Yau [Y] and S. Kolodziej [K 1], [K 2] that the equation $(\star)_t$ admits a unique continuous solution. (Observe that for $t \in ]0,1]$, $\omega_t$ is a K\"ahler form).
 
 Our aim here is to understand what happens when $t \to 0^+$, motivated by recent geometrical developpments  [ST], [KT]. When $n = m$, the cohomology class ${\om_0}$ is big and semi-ample and this problem has been adressed by several authors recently (see [CN], [EGZ], [TZ], [To]).
 
 We focus here on the case $m < n$. This situation is motivated by the study of the K\"ahler-Ricci flow on manifolds $X$ of intermediate  Kodaira dimension $1 \leq kod(X) \leq n - 1$. When $n = 2$ this has been studied by J.Song and G.Tian [ST].
 
 In a  very  recent and interesting paper [KT],  S. Kolodziej and G. Tian were able to show, under a technical geometric assumption on the fibration $\pi $, that 
 the solutions $(\vf_t)$ are uniformly bounded on $X$ when $t \searrow 0^+$.
 
 The purpose of this note is to (re)prove this result without any technical assumption and with a different method: we actually follow the strategy introduced by S. Kolodziej in [K] and further developped in [EGZ], [BGZ].
 
\vskip.2cm
 \noindent {\bf THEOREM.}
{\it There exists a uniform constant $M = M (\pi,\Vert F \Vert_p) > 0$ such that the solutions to  the 
 Monge-Amp\`ere equations $(\star)_t$ satisfy }
  $$  \Vert  \vf_t  \Vert_{L^{\infty}(X)} \ \leq \ M, \ \forall t \in ]0,1].$$


It follows from our result that Theorems 1 and 2 in [KT] hold without any 
technical assumption on the fibration (see condition 0.2 in [KT]).
\vskip.1cm

  This result has been announced by   J-P. Demailly and N. Pali [DP].
 
 \section{Proof of the theorem}
 
 \subsection{Preliminary remarks}
 
 \noindent {\bf Uniform control of $c_t$.} Observe that $\om_0^k = 0$ for $m < k \leq n$, hence for all $t \in ]0, 1]$, 
 $$ \omega_t^n = \sum_{k = 1}^m {n \choose k} t^{n - k} {\om_0}^k \we {\om_X}^{n - k}.$$
 Note that $]0,1] \ni  t \lto t^{m - n} \om_t^n$ is increasing (hence decreases as $t \searrow 0^+$) and satisfies for $t \in ]0,1]$ 
 
 \begin{displaymath}
 {n \choose m} \qquad \frac{\om_0^{m} \we \om_X^{n - m}}{\int_X \om_0^{m} \we \om_X^{n - m}} \leq \frac{\om_t^n}{t^{n-m} \int_X \om_0^m \we \om_X^{n - m}} \leq \frac{\om_1^n}{\int_X \om_0^m  \we \om_X^{n - m}}. \leqno (1)
 \end{displaymath}
 In particular $t \lto c_t$ is increasing in $ t \in ]0 , 1]$ and 
 $$ 0 <  {n \choose m} \frac{\int_X \om_0^m  \we \om_X^{n - m}}{\int_X F \om_X^n} = : c_0 \ \ \leq \ \ c_t \leq \ \ c_1.$$
\noindent {\bf Uniform control of densities.}
 Let $J_{\pi}$ denote the (modulus square) of the Jacobian of the mapping $\pi$, defined through   
 $$ \om_0^m \we \om_X^{n - m} = J_{\pi} \om_X^n.$$
 Let us rewrite the equation $(\star)_t$ as follows
 $$
 (\omega_t + dd^c \vf_t)^n  = f_t \omega_t^n,
 $$
 where for $t \in ]0,1]$
 $$0 \leq f_t := c_t t^{n - m} F  \frac{\om_X^n}{\om_t^n} \leq c_1 \frac {F}{J_{\pi}}.$$

Observe that 
$$
\int_X f_t \om_t^n = c_t t^{n - m} \int_X F \om_t^n 
=\int_X \om_t^n=: Vol_{\om_t}(X),
$$
hence $(f_t)$ is uniformly bounded in $L^1 (\om_t/V_t)$,
$V_t:=Vol_{\om_t}(X)$.
We actually need a slightly stronger information.
  
 \begin{lem} There exists $p' > 1$ and a constant $C = C (\pi,\Vert F\Vert_{L^p (X)}) > 0$ such that for all $t \in ]0,1]$
 $$ \int_X f_t^{p'} \om_t^n \leq C  \, Vol_{\om_t} (X).$$
 \end{lem}
 \noindent{\bf Proof of the lemma.} Set $V_t := Vol_{\om_t} = \int_X \om_t^n$ and observe that
 $$ 0 \leq f_t \frac{\om_t^n}{ V_t} \leq c_1 F  
 \frac{\om_X^n}{\int_X \om_0^m \we \om_X^{n - m}} = C_2 F  \om_X^n,$$
where $C_2 := c_1 \int_X J_{\pi} \om_X^n$.
 
 This shows that the densities $f_t$ are uniformly in $L^1$ w.r.t. the normalized volume fomrs $\om_t^n \slash V_t$.
 
 Since $J_{\pi}$ is locally given as the square of the modulus of a holomorphic function which does not vanish identically,  there  exists $\al \in ]0, 1[$ such that $ J_{\pi}^{- \al} \in L^1 (X)$.  Fix $\be \in ]0, \al[$ satisfying the condition $\be \slash p + \be \slash \al = 1$.  It follows from H\"older's inequality that
 $$ 
\int_X f_t^{\be} \om_X^n \leq \Bigl(\int_X F^p \om_X^n \Bigr)^{\be \slash p} \Bigl(\int_X J_{\pi}^{- \al}\omega_X^n \Bigr)^{\be \slash \al}.
$$
 Setting $\ep := \be \slash q$ and using H\"older's inequality again , we obtain
 $$\int_X f_t^{1 + \ep} \frac{\om_t^n}{V_t} \leq C_2
 \int_X f_t^{\ep} F \om_X^n.$$
 Now applying again H\"older inequality we get
 $$\int_X f_t^{1 + \ep} \frac{\om_t^n}{V_t} \leq C_2
 \Bigl(\int_X f_t^{\beta} \omega_X \Bigr)^{1 \slash q} \Vert F\Vert_{L^p (X)}.$$
 Therefore denoting by $p' := 1 + \ep$, we have the following uniform estimate
 $$\int_X f_t^{p'} \frac{\om_t^n}{V_t } \leq C (\pi,\Vert F\Vert_{L^p (X)}), \forall t \in ]0,1],$$
 where 
 $$C (\pi,\Vert F\Vert_{L^p (X)}) := C_2 \Bigl(\int_X J_\pi^{- \al}\omega_X^n \Bigr)^{\be \slash \al q}  \Vert F\Vert_{L^p (X)}^{1 + \beta \slash q }.$$
\hfill \fin
 
 \subsection{Uniform domination by capacity}

 We now show that the measure $ \mu_t := f_t \omega_t^n \slash Vol_{\om_t}$  
 are uniformly strongly dominated  by the normalized capacity  
  $\mathrm {Cap}_{\om_t}  \slash Vol_{\om_t} (X).$
 It actually follows from a carefull reading of  the no  
 parameter  proof given in [EGZ], [BGZ].  
 \begin{lem} There exists a constant $C_0 = C_0 (\pi,\Vert F\Vert_{L^p 
 (\om_X^n)}) > 0$ such that for any compact set $K \sub X$ and $t 
 \in ]0,1]$,
 $$ \mu_t (K) \leq C_0^n 
 \left(\frac{\mathrm {Cap}_{\om_t} (K)}{Vol_{\om_t}(X)}\right)^2.$$
 \end{lem}
 \demo Fix a compact set $K \sub X$.  Set $V_t := Vol_{\om_t} (X)$. H\"older's inequality yields 
 $$ \mu_t (K) \leq \left(\int_X f_t^{p'} \frac{\om_t^n}{V_t}\right)^{1 \slash p'} \left(\int_K \frac{\om_t^n}{V_t }\right)^{1 \slash q'}.$$
 It remains to dominate uniformly the normalized volume forms $ \om_t^n \slash V_t$  by the normalized capacities $\mathrm {Cap}_{\om_t} \slash V_t$. Fix $\si > 0$ and observe that for any $t \in ]0,1]$,
 $$ \int_K \frac{\om_t^n}{V_t} \leq \int_X e^{- \si ( 
 V_{K,\om_t} - \max_X V_{K,\om_t})} \frac{\om_t^n}{V_t} 
 T_{\om_t} (K)^{\si},$$
 where 
$$V_{K,\om_t} := \sup \{ \psi \in PSH (X,\om_t) ; \psi \leq 0, \ \mathrm {on} \ K\}$$
 is the $\om_t-$extremal function of $K$ and  
 $T_{\om_t} (K) := \exp (- \sup_X V_{K,\om_t})$ is the associated $\om_t-$capacity of $K$ (see [GZ 1] for their properties).   
 
 Observe that $\om_t^n \slash V_t \leq c_1 \om_1^n$ and $\om_t 
 \leq \om_1$, hence the family of functions $V_{K,\om_t} - \max_X 
 V_{K,\om_t}$ is a normalized family of $\om_1-$psh functions. Thus there exists $\si > 0$ which depends only on $(X,\om_1)$ and 
 a constant $ B = B (\si,X,\om_1)$ such that ([Z])
 $$ \int_X e^{- \si ( V_{K,\om_t} - \max_X V_{K,\om_t})} 
 \frac{\om_t^n}{V_t} \leq B, \forall t \in 
 ]0,1].$$
 The Alexander-Taylor comparison theorem (see Theorem 7.1 in [GZ 1]) now yields for a constant $C_3 = C_3 (\pi, \Vert F\Vert_{L^p(X)}) $
 $$ \mu_t (K) \leq  C_3 \exp \left[- \si \left(\frac{V_t}{\mathrm{Cap}_{\om_t} (K)}\right)^{1 \slash n}\right], \forall t \in ]0,1].$$
 We infer that  there is a constant $C_4 = C_4 (\pi, \Vert F\Vert_{L^p(X)})$ such that

 $$ 
\mu_t (K) \leq C_4 \left(\frac{\mathrm{Cap}_{\om_t} (K)}{V_t}\right)^2, \forall t \in ]0,1]. 
\leqno (2)
$$

\subsection{Uniform normalization}

 The comparison principle (see [K] ,[EGZ]) yields for any $s > 0$ and $\tau \in [0,1]$
 $$ \tau^n \frac{\mathrm{Cap}_{\om_t} (\{\vf_t \leq - s - \tau\})}{V_t}
 \leq \int_{\{\vf_t \leq - s \}} \frac{(\omega_t + dd^c \vf_t)^n}{V_t}.$$
 It is now an exercise to derive from this inequality  an a priori  $L^{\infty}-$estimate,
 $$
\Vert \vf_t \Vert_{L^{\infty} (X)} \leq C_5 + \ {s}_0 (\om_t),
$$
 where $ {s}_0 (\om_t) $  (see [EGZ],[BGZ]) is  the smallest number $s > 0$ satisfying the condition
 $ e^n C_0^n  \mathrm{Cap}_{\om_t}(\{\psi \leq - s \})/V_t \leq 1$ for all $\psi \in PSH (X,\om_t)$ such that 
$\sup_X \psi = 0$.
 Recall from ([GZ 1], Prop. 3.6) that
 $$ 
\frac{\mathrm{Cap}_{\om_t} (\{\psi \leq - s - \tau\})}{V_t}
 \leq \frac{1}{s} \left(\int_X (- \psi) \frac{\om_t^n}{V_t} + n \right).
 $$
 Since $\frac{\om_t^n}{V_t} \leq C_1 \om_1^n$, it follows that
 $$ 
\frac{\mathrm{Cap}_{\om_t} (\{\psi \leq - s - \tau\})}{V_t}
 \leq \frac{1}{s} \left(C_1 \int_X (- \psi) \om_1^n + n \right).
$$
 Since $\psi$ is $\om_1-$psh and normalized, we know that there is a constant $A = A (X,\om_1) > 0$ such that $C_1 \int_X (- \psi) {\om_1^n} \leq A $ for any such $\psi$. Therefore  $s_0 (\om_t) \leq  s_0 := e^n C_0^n (A  + n)$ for any $t \in ]0,1]$. 
 Finally we obtain the required uniform estimate for all $t \in ]0,1]$.

 \end{document}